\documentclass{article}
\usepackage{amsmath,amscd,amsfonts,amssymb,amsthm,newlfont,exscale,color,mathrsfs, multicol}

\title{Ideals Generated by Quadratic Polynomials}
\author{Tigran Ananyan and Melvin Hochster\footnote{The second author was partially supported by a grant from the National Science Foundation (DMS--0901145).}}
\date{\today}

\usepackage{enumerate}
\newtheorem{conjecture}{Conjecture}
\newtheorem{CastelnuovoMumford}[conjecture]{Conjecture}
\newtheorem{MainTh}{Main Theorem}
\newtheorem{MainTh2}[MainTh]{Main Theorem}
\newtheorem{HeightDoesNotRise}{Lemma}
\newtheorem{RegSequenceLemma}[HeightDoesNotRise]{Lemma}
\newtheorem{RegSeqLemma}[HeightDoesNotRise]{Lemma}
\newtheorem{LinLemma}[HeightDoesNotRise]{Lemma}
\newtheorem*{keylemma}{Key Lemma}
\newcommand\blank{\underline{\hbox{\ \ }}}
\newcommand\exn{coefficient\ }
\newcommand\pf{\noindent \textit{Proof.}\ }
\newcommand\inc{\subseteq}
\newcommand\uu{\underline{u}}
\newcommand\uv{\underline{v}}
\newcommand\uw{\underline{w}}
\newcommand\uz{\underline{z}}
\newcommand\std{standard form}
\newcommand\vect[2]{#1_1,\,\ldots,\, #1_{#2}}
\newcommand\ve[3]{#1_{#2}, \, \ldots, \, #1_{#3}}
\begin{document}
\maketitle

\begin{abstract} Let $R$ be a polynomial ring in $N$ variables over an arbitrary field $K$ and let $I$ be
an ideal of $R$ generated by $n$ polynomials of degree at most 2.  We show that there is a bound
on the projective dimension of $R/I$ that depends only on $n$, and not on $N$.  The proof depends
on showing that if $K$ is infinite and $n$ is a positive integer, there exists a positive integer $C(n)$,
independent of $N$, such that any $n$  forms of degree at most 2 in $R$ are contained in a subring
of $R$ generated over $K$ by at most $t \leq C(n)$ forms $G_1, \, \ldots, \, G_t$ of degree 1 or 2
such that $G_1, \, \ldots, \, G_t$ is a regular sequence in $R$. $C(n)$ is asymptotic to $2n^{2n}$.
\end{abstract}

\section{Introduction}

Throughout this paper let $R$ denote a polynomial ring over an arbitrary field $K$:  say
$R = K[x_1, \ldots , x_N]$. We will denote the projective dimension of the
$R$-module  $M$ over $R$ by $\textrm{pd}(M)$. The following conjecture was posed
by M.~Stillman in \cite{PS}.

\begin{conjecture}\label{mainconj}
There is an upper bound, independent of $N$, on $\textrm{pd}(R/I)$, where $I$ is any ideal of $R$
generated by $n$ homogeneous polynomials of given degrees $d_1,\,\ldots,\, d_n$.
\end{conjecture}

One motivation for proving this conjecture comes from its equivalence to the following open question
(the proof of equivalence due to Caviglia is given in \cite{E1}):

\begin{CastelnuovoMumford}
There is a bound on the Castelnuovo-Mumford regularity of an ideal in  polynomial ring that depends
only on the number of its minimal generators and the degrees of those generators.
\end{CastelnuovoMumford}

Another motivation for Conjecture \ref{mainconj} comes from certain results about ideals generated by
three homogeneous polynomials. A construction of Burch (\cite{BU}) in the local case, extended to the
global case by Kohn (\cite{Ko}), shows that there exist three-generated ideals of arbitrarily large projective
dimension. However in the polynomial ring case, as the projective dimension grows, the degrees of the
generators are growing as well, thus motivating Conjecture \ref{mainconj} (see \cite{E1} for analysis of
growth of degrees in Burch's construction).\\

Conjecture \ref{mainconj} clearly holds when $n \leq 2$ or when all of the $d_i = 1$,
but even the simplest next case of $n = 3$ and $d_1 = d_2 = d_3 = 2$ requires non-trivial arguments.
D.~Eisenbud and C.~Huneke proved that if $I$ is generated by three quadratic forms, then
$\textrm{pd}(R/J) \leq 4$. B.~Engheta showed the existence of a bound on the projective dimension for
three cubic forms (see \cite{E2}, \cite{E3}); he has shown that $\textrm{pd}(R/I)\leq 36$ (even though in
all known examples of three cubics $\textrm{pd}(R/I)\leq 5$). In general, the projective dimension can
grow relatively fast as one increases the number of generators and the degrees (see \cite{BS}, \cite{CK},
and \cite{Mc} for specific examples). In this paper the authors prove Conjecture \ref{mainconj} for the case
when all of the $d_i$ are at most $2$ and $n$ is arbitrary.\\

Since a base change of the field $K$ to a larger field does not affect the projective dimension, in the
remainder of this paper we often pass to the case  where $K$ is infinite.\\

\section {The Main Results}  We use recursion on $h$ to define a function $B(m,n,h)$ of three nonnegative integers, with $h \leq n$, as follows.
Let $B(m,n,0) = m +mn=m(n+1)$.  If  $h \geq 1$,  let $$B(m,n,h)= (m+h)(n^3+n^2 +n +1) + h(n+1) + B\bigl((m+h)n^2, n, h-1\bigr).$$

We also define $C(s)$ to be the largest value of $B(m,n,h)+h$ for $m+n = s$ and $0 \leq h \leq n-1$.
 Let $C_0(s)$ be the largest value of $B(0,n,h)+m+h$ for $m+n=s$ and $0 \leq h \leq n-1$. It is shown in \S\ref{size} that both $C_0(s)$ and
 $C(s)$ are asymptotic to $2s^{2s}$.\\

Our goal is to prove the following result.

\begin{MainTh}  Let $R=K[x_1,\, \ldots,\, x_N]$ be a polynomial ring over a field $K$.
Let  $F_1,\, \ldots,\, F_{m+n}$ be forms of degree at most 2 generating ideal $I$, and suppose that the linear forms among them
span a $K$-vector space of dimension $m$.  Let  $h$ denote the height of the ideal  generated by
the images of the quadratic $F_j$ in the ring obtained from $R$ by killing the ideal generated by all of the $F_i$
that are linear.

If $K$ is infinite, after a linear change of variables, there are at most $b \leq B(m,n,h)$   variables $\vect y b$
and at most $c \leq h$ quadratic forms  $\vect G c$ in $I$ such that
\begin{enumerate}[(1)]
\item  The elements $\vect y b,\,\vect G c$ form a regular sequence in $R$.
\item  The $F_i$ are in the polynomial ring  $K[\vect y b,\,\vect G c]$.
\item  The quadratic forms $\vect G c$ are also a regular sequence modulo those generators  of $I$  that
involve only the variables $\vect y b$.
\end{enumerate}
(Condition (3) is automatic: see Lemma \ref{regseq} below).

Hence, if $K$ is infinite,  the elements $F_1,\, \ldots,\, F_{m+n}$ are contained in the $K$-subalgebra of $R$
generated by a regular sequence of at most $B(m,n,h) + h$ linear and quadratic forms.  It also follows that, if
$K$ is infinite,  the elements $F_1,\, \ldots,\, F_{m+n}$ are contained in the $K$-subalgebra of $R$ generated
by a regular sequence of at most $C(m+n)$ forms of degree at most 2.

Consequently, for every field $K$,  the projective dimension of $R/I$  is at most  $B(0,n,h) +m+h$,
and is also at most $C_0(m+n)$.

\end{MainTh}

Note that the last statement follows from the first because the projective dimension drops by $m$ when we kill $m$ variables
in $I$ and pass to the polynomial ring in the remaining variables, replacing the original $F_j$ by the images of the quadratic
$F_j$ in the smaller polynomial ring.

Note also that if $h=n$ then there is a much better result, in that $\vect F {m+n}$ is already a regular sequence.

Observe also that, as mentioned earlier, the  projective dimension does not change when we enlarge the base field.\\

Given polynomials of degree
at most 2 that need not be forms, each is the sum of at most one quadratic form, one linear form, and a scalar
that is already in $K$.    Hence:

\begin{MainTh2} Let $R=K[x_1,\, \ldots,\, x_N]$ be a polynomial ring over a field $K$. Let  $F_1,\, \ldots,\, F_s$
be polynomials of degree at most 2, and let $I = (F_1,\, \ldots, \,F_s)$. If $K$ is infinite, the  elements
$F_1,\, \ldots,\, F_s$ are contained in the $K$-subalgebra of $R$ generated by a regular sequence of at
most $C(2s)$ forms of degree at most 2.

Consequently, for any field $K$, the projective dimension of $R/I$  is at most $C(2s)$. $\square$
\end{MainTh2}

Finally, observe that these theorems hold even if we allow the ambient polynomial ring to have a set of variables of arbitrary
infinite cardinality, because any given finite set of polynomials will only involve finitely many of the variables.

\section{Preliminary results}

To prove these theorems, we will need the following four lemmas.

\begin{HeightDoesNotRise}
Let $I$ be a homogeneous ideal in the polynomial ring $R$. If we take the image of $I$ after killing some of
the variables,  its height does not increase.  Much more generally, let $R$ be a locally equidimensional catenary
Noetherian ring, let $x$ be any element not in any minimal prime of $R$, and let $I$ be any ideal such that for
every minimal prime $P$ of $I$, $P + xR \not= R$.  Then the height of $I(R/xR)$ is at most the height of $I$.
\end{HeightDoesNotRise}

\textit{Proof.} The first statement follows from the second by induction on the number of variables killed.  Note
that when $I$ is proper homogeneous, so are all of its minimal primes, and this is preserved both by adding
$xR$ for a variable $x$ and by killing $xR$. Without loss of generality, we may replace the ideal $I$ by a minimal
prime $P$ of the same height. Since $R$ is locally equidimensional and catenary, the proper ideal $P+xR$ has
height equal to either $\textrm{ht}(P)$ or  $\textrm{ht}(P)+1$, depending on whether $x \in P$ or not:  this can be
checked after localizing at any minimal prime of $P+xR$. After $xR$ is killed,
the height drops by $1$, because $x$ is not in any minimal prime of $R$. $\square$

\begin{RegSequenceLemma}\label{regseq}
If  $F_{r+1}, \, \ldots, \, F_{r+s}$ form a regular sequence modulo a set of variables containing  the variables occurring in
$F_1, \ldots, F_r$,  then they form a regular sequence modulo the ideal $(F_1, \ldots, F_r)$.
\end{RegSequenceLemma}

\pf The variables killed generate a ring $A$,  and we may form the quotient ring  $B = A/(F_1, \ldots, F_r)A$. The polynomial
ring $C$ in the rest of the variables over $B$ is flat over $B$. It suffices to prove the statement after replacing $B,\, C$
by their localizations at their homogeneous maximal ideals.  But then the result follows from \cite{Mat}, Ch. 8, (20.F): if
the image of a sequence in the closed fiber of a flat local extension is regular, so is the original sequence. $\square$\\

A regular sequence of forms in a polynomial ring over a field $K$ can be extended to a homogeneous system of
parameters.  The original polynomial ring is module-finite and free over the polynomial ring generated over $K$ by the
homogeneous system of parameters. Hence:

\begin{RegSeqLemma}\label{regsequ} Let $F_1,\,\ldots,\,F_t$ be a regular sequence of forms in a polynomial ring $R$
over a field $K$. Then $R$ is free, hence, faithfully flat over $A = K[F_1,\,\ldots,\,F_t]$.  Hence, for any ideal $J$ of $R$
whose generators lie in $A$,  $\textrm{pd}(R/J) \leq t$.
$\square$
\end{RegSeqLemma}

\begin{LinLemma}\label{lin} Let $K$ be any field, let $\vect y r,\,\vect z s$ be $r+s$ indeterminates over $K$,  let
$\vect \alpha  n \in K[\vect y r]$,  and let $\vect \beta  n \in K[\vect z s]$.  Then we can map the polynomial ring
$K[\vect Tn]$ onto $K[\vect \alpha n]$ as a $K$-algebra so that each $T_i \mapsto \alpha_i$.  Call the kernel  $P$.
Likewise, we can map the polynomial ring  $K[\vect Tn]$ onto $K[\vect \beta n]$ as a $K$-algebra so that each
$T_i \mapsto \beta_i$.  Suppose that the kernel of this map is also  $P$.  Finally, we can map the polynomial ring
$K[\vect Tn]$ onto $$K[\alpha_1+\beta_1, \, \ldots, \, \alpha_n+\beta_n] \inc K[\vect yr,\,\vect zs]$$
as a $K$-algebra so that each $T_i \mapsto \alpha_i+\beta_i$.  Suppose that the kernel of this map contains  $P$.
Assume also that $P$ is homogeneous, which is automatic if the $\alpha_i$ (or the $\beta_i$) are homogeneous
of the same degree over $K$.  Then $P$ is generated by linear forms over $K$.
\end{LinLemma}

\pf Let $L$ denote an algebraic closure of $K$.  Since $$0 \to P \to K[\vect Tn] \to K[\vect \alpha n] \to 0$$
is exact and $K[\vect \alpha n] \to K[\vect y r]$ is injective,  the same holds after we apply $L\otimes_K \blank$. Similar
remarks apply when we consider the $K$-algebra spanned by the $\beta_i$.  It follows that  $L\otimes_K P$, which may
be identified with expansion of $P$ to $L[\vect y r]$, is the kernel of the map of $L[\vect T n]$ onto $L[\vect \alpha n]$,
and taken together with the corresponding facts for the $\beta_i$ and the  $\alpha_i+\beta_i$, this shows that our hypotheses
are preserved after base change to $L$.  If we know the result for $L$, we know that  $L\otimes_KP$ is generated as an
ideal by its degree 1 part, $[L\otimes_KP]_1 = L\otimes_K[P]_1$.  Then   $[P]_1K[\vect T n]$ expands to $PL[\vect Tn]$.
Since $L[\vect T n]$ is faithfully flat over $K[\vect T n]$, the contractions  $[P]_1K[\vect T n]$ and $P$, respectively, are
equal. That is,  $P$ is generated by $[P]_1$.

Hence, we may assume without loss of generality that  $K$ is algebraically closed.  Consider the algebraic set $V$ defined by
$P$ in affine $n$-space over $K$.  $V$ is closed under scalar multiplication because $P$ is homogeneous.  We claim that
$V$ is closed under addition.  For suppose that $(\vect c n),\, (\vect d n) \in V$.  Then there is a $K$-homomorphism
$K[\vect \alpha n] \to K$ sending each  $\alpha_i \to c_i$ and, similarly, a $K$-homomorphism $K[\vect \beta n] \to K$  sending
every $\beta_i \to d_i$.  This yields a $K$-homomorphism $\theta:K[\vect \alpha n] \otimes_K K[\vect \beta n] \to K$ such that
for every $i$, $\alpha_i \mapsto c_i$ and $\beta_i \mapsto d_i$.  We may tensor the injections  $K[\vect \alpha n]  \inc K[\vect y r]$
and  $K[\vect \beta n] \inc K[\vect zs]$ over $K$  to obtain an injection  $$K[\vect \alpha n] \otimes_K K[\vect \beta n] \inc
K[\vect yr, , \, \vect zs]$$ that enables us to identify  $K[\vect \alpha n] \otimes_K K[\vect \beta n]$ with its image
 $K[\vect \alpha n, \vect \beta n] \inc K[\vect yr, \, \vect zs]$.  Hence, we may think of  $\theta$ as a $K$-homomorphism
 $K[\vect \alpha n, \vect \beta n] \to K$  such that  for every $i$, $1 \leq i \leq n$,   $\alpha_i + \beta_i \mapsto c_i+ d_i$.   Since every polynomial
 $H \in P$ vanishes on $(\alpha_1 + \beta_1,\, \ldots,\, \alpha_n + \beta_n)$,  we have that $H(c_1+d_1, \, \ldots, \, c_n+d_n)= 0$
 as well.  This completes the proof that $V$ is closed under addition, and so $V$ is a vector subspace of $K^n$.  But then its defining
 prime ideal $P$ is generated by linear forms. $\square$

\section {The basic step: putting the $F_i$ and variables in \std.}

We are going to make iterative of use the construction of putting the $F_i$ and the variables in \std\ described just below. \\

We shall say that the $F_i$ and the variables $x_j$  are in {\it \std} if the following conditions hold:

\begin{enumerate}[(1)]

\item  $F_{n+i} = x_i$ for $1\leq i \leq m$ and the remaining $F_i$ are quadratic or 0.

\item No $F_i$ for $1 \leq i \leq n$  has a monomial term in $K[\vect x m]$.

\item $F_1, \, \ldots F_h,\, \vect x m, \, x_{m+h+1}, \, \ldots, \, x_N$  is a regular sequence.

\item There is a nonnegative integer $r \leq (m+h)n$ such that when the $F_i$ are written as polynomials in the variables $\vect x {h+m}$
with coefficients in the ring $K[x_{h+m+1},\, \ldots, \, x_N]$,  the coefficients occurring for the variables $\vect x {h+m}$ are in
the $K$-span of $x_{h+m+1},\, \ldots, \, x_{h+m+r}$.

\item There is a nonnegative integer $s \leq (m+h)n^2$ such that  when the $F_i$ are written as polynomials in the variables $\vect x {h+m+r}$
with coefficients in the ring $K[x_{h+m+r+1},\, \ldots, \, x_N]$,  the coefficients occurring for the variables $x_{h+m+1},\,\ldots, x_{h+m+r}$ are in
the $K$-span of $x_{h+m+r+1},\, \ldots, \, x_{h+m+r+s}$.

\item  We call $\ve x {m+1} {m+h}$ {\it front variables}. We shall also denote them $\vect u h$.
Let $f_i$ denote the image of $F_i$ under the $K$-homomorphism
$\pi:R \to K[\vect u h]$ that kills all the variables except the front
variables while fixing the front variables (note that $\pi$ kills $F_i$ for $i > n$). Then
there is an integer $d$, $h \leq d \leq n$, such that $\vect f d$  are linearly independent over $K$  while $f_i = 0$
for $i > d$. We call $\vect f n$ the {\it front polynomials}.

\end{enumerate}

Under these conditions we call $\vect xm$ {\it leading variables},  and, as already mentioned, we call $\ve x {m+1} {m+h}$ the front variables,
and denote them $\vect u h$.    We call $\ve x {h+m+1} {h+m+r}$
the {\it primary \exn   variables} and use the alternative notation $\uv=\vect v r$ for them.  We call $$\ve x {h+m+r+1} {h+m+r+s}$$
the {\it secondary \exn   variables} and we write $\uw =\vect w s$ as an alternative notation for them.  We refer to $\ve x {h+m+r+1} N$ as
the {\it tail variables}.   We use the alternative notation $\uz = \vect z {N-(h+m+r)}$ for the tail variables.
It is very important to note that the tail variables include the secondary \exn   variables $\uw$.

Let $\tau$ denote the $K$-homomorphism from $R$ to the
polynomial ring in the tail variables that fixes the tail variables and kills the others.  We write $g_i = \tau(F_i)$.
Again, $\tau$ kills $F_i$ for $i > n$. We call $\vect gn$  the {\it tail polynomials}.

We can map the polynomial ring $K[\vect T n]$ onto $K[\vect fn]$ by $T_i \mapsto f_i$.  The kernel is a prime ideal $P$
in $K[\vect Tn]$, which we call the ideal of {\it front relations}.  \\

We want to see that these conditions can be achieved.  Note, however, that \std\ is far from unique.\\

\noindent{\bf Achieving standard form.}\\

We want to show that the sequence can be put in standard form without changing the ring it generates.

If the $F_i$ are linearly dependent over $K$, one or more $F_i$ may be replaced by 0,  and we assume that
any 0 elements occur in the final spots among the first $n$ terms of the sequence. Hence,
we may assume that the linear forms that occur are linearly independent.
We may assume that the $F_i$ are numbered so that any linear forms occur
as a final segment, and we may assume these forms are also an
initial segment of the variables. Thus, we have  $x_i = F_{n+i}$, $1 \leq i \leq m$.

We may now subtract from each $F_i$, $1 \leq i \leq n$, the sum of all terms that occur that involve only
$\vect x m$.

Since $K$ is infinite, after replacing  $F_1,\,\ldots,\, F_h$  by suitably general linear combinations of the quadratic forms
$F_1,\,\ldots,\, F_n$, we may assume that the images of $F_1,\, \ldots, \,F_h$ form a maximal regular sequence in
$(\vect F n)R/(\vect x m)R$.  Note that for $i \leq n$,  every $F_i$ will be quadratic or else 0.  Since regular sequences of forms
are permutable, we have that
$F_1,\,\ldots,\, F_h,\,\vect x m$ is a regular sequence, and since $K$ is infinite
we may extend $F_1,\, \ldots, \,F_h, \vect x m$ to a homogeneous system of
parameters using linear forms,  which we may assume are the variables $x_{m+h+1}, \, \ldots,\, x_N$.
That is, $F_1,\, \ldots, \,F_h,\, \vect x m,\, x_{h+m+1}, \, \ldots,\, x_N$  form a homogeneous system of parameters
for $K[x_1, \ldots, x_N]$.  The images of $\vect F h$  when we kill the variables in this regular sequence
are precisely the elements $\vect f h \in K[\vect u h]$.  Hence, the quadratic forms $\vect fh$ are a homogeneous
system of parameters for $K[\vect uh]$.  In particular, they are linearly independent over $K$. The remaining $F_i$
may be permuted so that $\vect f d$  is a $K$-vector space basis for the $K$-span of the $f_i$.  For $i > d$ we
may subtract a $K$-linear combination of $\vect F d$  from $F_i$ to arrange that $f_i$ be zero.

When $\vect F n$ are written as  polynomials in the leading and front variables $\vect xm,\,\vect uh$ (with coefficients in the polynomial ring
in the remaining variables) there are most  $m+h$ linear terms:  the coefficients are linear forms in $\ve  x {m+h+1} N$.  All these coefficients from all of the
$F_i$ span  a $K$-vector space of dimension $r \leq n(m+h)$.  We make a linear change of the variables $\ve  x {m+h+1} N$
so that the variables $\ve x {m+h+1}  {m+h+r}$ span this space.  These are the primary \exn   variables, which are also denoted
$\vect v r$.

Next, when $\vect F n$ are written as  polynomials in the leading, front, and primary \exn   variables  $\vect xm,\,\vect uh,\,\vect vr$
(with coefficients in the polynomial ring in the remaining variables) there are most  $n(m+h)$ terms that are linear involving one of the $v_j$:
the coefficients are linear forms in $\ve  x {m+h+r+1} N$.  All these coefficients from all of the
$F_i$ span  a $K$-vector space of dimension $s \leq n^2(m+h)$.  We make a linear change of the variables $\ve  x {m+h+1} N$
so that the variables $\ve x {m+h+r+1}  {m+h+r+s}$ span this space.  These are the secondary \exn   variables, which we also
denote $\vect w s$.

Finally, we refer to $\ve x {m+h+r+1} N$ as the {\it tail variables}.  This set is the complement in $\{\vect x N\}$ of the union of
leading, front, and primary \exn   variables.  We emphasize again that the tail variables include the secondary \exn   variables.

All of the conditions for standard form are now satisfied.  We remark that once we have standard form, it is unaffected by
permuting  $F_{h+1},\,\ldots,\,F_d$, or by multiplying by an invertible matrix over $K$, thereby replacing them with linear
combinations that have the same $K$-span.

We show in the next section that if this procedure is then carried through a second time using the forms
consisting of all leading, front, primary and secondary \exn   variables (these are the new leading variables) and the tails, thus producing
a second set of tails, then one of two things happens:
\begin{enumerate}[(1)]

\item The second set of tails generates an ideal of height at most $h-1$ modulo the ideal generated by the new leading variables. (It may
also happen that the set of original tails generates an ideal of height at most $h-1$ modulo the ideal generated by the original leading variables,
which is an easier case.)
\item The leading, front, primary and secondary \exn   variables and the nonzero elements in first set of tails form
a regular sequence.

\end{enumerate}

Either condition yields an estimate of what is needed for $B(m,n,h)$.  The first condition requires a larger value and leads to the
recursive definition given earlier.  Details are given in the next section.

\section{Key behavior of \std\ and the proof of the main theorems}

Placing $\vect F {m+n}$ and the variables in \std\ puts surprising constraints on the forms.  In particular, parts (c), and (d)
of the Key Lemma below play a central role in the proof of the main results.

\begin{keylemma}\label{key} Let $\vect F {m+n}$ consist of quadratic and linear forms in the polynomial ring
$R =K[\vect x N]$ which, together with the variables  $\vect xN$ are assumed to be in \std. Let all notation and terminology, including $m,\,n,\,h,\, r,\,s,\,d,$
 $f_i, \, g_i, \, u_i,\,  v_i,\, w_i,\, P,\, \pi,$ and $\tau$, be as in the preceding section. Let $I = (\vect F {m+n})R$.  Then:

\begin{enumerate}[(a)]

\item The elements $\vect fh$ form a homogenous system of parameters for the ring $K[\vect uh]$, and since
$$K[\vect fh] \inc K[\vect fd] = K[\vect fn] \inc K[\vect uh]$$ the domain $K[\vect fn]$ has Krull dimension $h$,
so that $P$ has height  $n-h$.

\item For $1 \leq i \leq n$, $F_i$ uniquely has the form  $f_i +e_i +g_i$  where
$$e_i \in K[\vect xm, \,\uu,\, \uv,\, \uw]$$  and is also
in the ideal of this ring generated by $(\vect xm,\,\uv)$.  Hence,
$$K[\vect F {m+n}] \inc K[\vect g n,\, \vect xm,\,\uu,\, \uv,\, \uw].$$

\item  If $H \in P$,  then  $H(\vect g n) =0$.   That is, the tail polynomials satisfy the front relations.

\item   If $i > d$,  the  $g_i = 0$,  i.e.,  $F_i = e_i$ is in the ideal generated by the leading and primary \exn   variables,
and in the polynomial ring generated by the leading variables, front variables and the primary and secondary \exn   variables.

\end{enumerate}
\end{keylemma}

\pf  The statement that $\vect f h$ is a homogeneous system of parameters was proved in the fourth paragraph of the
subsection on Achieving standard form, and the other assertions in (a) are immediate.\\

For (b),  $e_i$ is clearly the sum of all terms in $F_i$ not involving only front variables nor only tail variables.
From the condition (2) of the definition of standard form  there are no terms involving only
leading variables. The terms in $e_i$
that involve some $x_i$ or $u_i$ in degree 1 have coefficients in the span of the $v_i$ by the definition of
primary \exn   variables.
That is,  $$e_i = \sum_{j=1}^m L_jx_j + \sum_{k=1}^h L'_k u_k + e'_i$$  where the $L_j,\, L_j'$
are in the $K$-span of $\vect vr$ and the terms in $e'_i$ are quadratic in the $v_t$ or linear
in the $v_t$,  and those that are linear in the $v_t$ have coefficients in the $K$-span
of $\vect w s$ by the definition of secondary \exn   variables.\\

To prove (c),  note first that $H$ is a sum of homogeneous elements  of  $P$ and so may be assumed homogeneous
of degree $k$.  Since $-f_i \equiv e_i + g_i$ modulo $I$,  we have
$$0 = (-1)^kH(\vect fn)= H(\vect {-f}n) \equiv H(e_1  + g_1,\,\ldots,  \,  e_n  + g_n)\hbox{\ modulo\ } I,$$  so that
$$H(e_1 +g_1,\, \ldots, \, e_n  + g_n) \in I.$$
If  $H(\vect g n) \not=0$, we obtain a contradiction by showing that  $$\vect F h,\, H(e_1  + g_1,\,\ldots,  \,  e_n  + g_n)$$
is a regular sequence of length $h+1$ in $R/(\vect xm)$.  To this end, it suffices to show that we have a regular
sequence modulo $(\vect x m, \vect v r)$.  Killing the $v_j$ kills the $e_j$ by (b) above,  and so it suffices to show
that $$f_1+g_1, \, \ldots, \, f_h + g_h, \, H(\vect g n)$$ is a regular sequence in the polynomial ring in the front
and tail variables.  This is immediate from Lemma \ref{regseq}:  since regular sequences of forms are permutable,
it suffices to show that $$H(\vect g n),\, f_1+g_1, \, \ldots, \, f_h + g_h$$ is a regular sequence.  Since the first element is nonzero, it is
a nonzerodivisor involving only the tail variables, and modulo the tail variables the remaining terms become the
regular sequence $\vect f h$.\\

Finally, (d) follows from (c) and the fact that $f_i = \pi(F_i) =0$ for $i > d$:  we may take  $H = T_i$, and this gives
that $g_i =0$.  The remaining statements now follow from (b) above.
$\square$\\

\noindent\textit{Proof of Main Theorem 1.} If $h=0$, then $\vect F n$ are expressible as sums of multiples by linear
forms of $\vect xm$.  All the multipliers together span a vector space of dimension at most  $mn$, and so all of
$\vect F n$ can be expressed in terms of $m + mn$ variables, which include the variables $\vect xm$.\\

We now assume $h\geq1$ and use induction on $h$.  We may assume standard form. We consider two cases, (1) and (2).  Case (2)
has subcases (2a) and (2b). Note that Case (1) is, in a sense, subsumed in Case (2a).
Let $h'$ be the height of the ideal generated by the tail polynomials
modulo the secondary \exn   variables.  Of course, $h' \leq h$.\\

\noindent {\bf Case 1:  $h'\leq h-1$.}  By the induction hypothesis, we may work in the polynomial ring in the tail
variables $\uz$ with the sequence of tail polynomials and secondary \exn   variables $\uw$ to select at most $k \leq B(n^2(m+h), n, h-1\bigr)$
tail variables $z_{j_1},\,\ldots, z_{j_k}$, such that $K[\vect gn, \, \vect w s]$ is contained in the $K$-algebra generated by these $z_{j_\nu}$ and a sequence of
at most $h-1$ distinct quadratic forms  that, together with the $z_{j_\nu}$ form a regular sequence. When we include the leading, front, and primary \exn   variables (which are disjoint from the tail variables), we still have a regular sequence,   and by part (b) of the Key Lemma,  the algebra these generate
contains $K[\vect F {m+n}]$. In this case the constraint placed on $B(m,n,h)$  is that it be at least
$$m+h+n(m+h)  + B\bigl(n^2(m+h), n, h-1).$$ \\

\noindent {\bf Case 2:  $h' = h$.} Let $S$ be the polynomial ring over $K$ generated by the tail
variables. We work with sequence $\vect g n,\, \vect ws$ in $S$.  We may now put these in standard form.  We can do so without
changing $\vect g h$:  we claim these already form a regular sequence of length $h$ modulo the secondary \exn   variables. (To meet the
condition that an initial segment of the new front polynomials be independent over $K$ while the rest are zero, we may adjust
$F_{h+1},\, \ldots, \, F_d$ by multiplying by an invertible matrix over $K$: note the comment on p.~7.)  
The reason is that there is a homomorphism
of $K$-algebras $K[\vect f n] \to K[\vect g n]$  sending $f_i \mapsto g_i$  for all $i$.  This is well-defined precisely because of
part (d) of the Key Lemma:  the $g_i$  satisfy all the relations in $P$ on the $f_i$.  It is clearly surjective. it follows that
$K[\vect g n]$ is module-finite over $K[\vect g h]$.  But $K[\vect g n]$ maps onto the new ring of front polynomials,
which will have Krull dimension $h$.  Thus,  $K[\vect g h]$ must have dimension $h$, and no less. This implies that
$\vect g h$ is a regular sequence, even modulo the secondary \exn   variables.  Thus, the prime ideal of relations
on $\vect g n$ must be exactly $P$.  The same holds for the new ideal of front relations on the front polynomials of the $g_i$.
We let $\vect \alpha n $ denote the new front polynomials let  $\vect \beta n$ denote the new tail polynomials.
Let $h''$ denote the height of the ideal generated by the new tail polynomials in the ring generated over $K$ by the new
tail variables modulo the ideal generated by the new secondary \exn   variables. \\

\noindent {\bf Subcase (2a):  $h'' \leq h-1$.} This is very much like Case (1).  We have that $K[\vect F {m+n}]$
is contained in the polynomial ring generated over $K$  by the original leading, front, primary \exn   variables,
secondary \exn   variables, and the original tail polynomials.  But  $K[\vect ws, \, \vect g n]$
is contained in turn in the $K$ algebra generated by the secondary \exn   variables (which are also the
new leading variables), the new front, new primary \exn  , and new secondary \exn   variables, and the new tails.
We count variables as follows:

$$\begin{array}{ll}
\underline{\hbox{\rm Type of variable or term}}       &  \underline{\hbox{\rm  Cardinality}}\\
{}&{}\\
\hbox{\rm leading}                                                                    &                  m\\
\hbox{\rm front}                                                                         &                 h\\
\hbox{\rm primary \exn  }                                                 &           (m+h)n\\
\hbox{secondary \exn   = new leading}                     &           (m+h)n^2\\
\hbox{\rm new front}                                                                 &               h\\
\hbox{\rm new primary \exn  }                                         &           \bigl((m+h)n^2+h\bigr)n \\
\hbox{\rm new secondary \exn   and new tails}               &     \hbox{\rm covered by } B\bigl((m+h)n^2,n,h-1\bigr)\\
\end{array}$$
In this case the constraint placed on $B(m,n,h)$ is that it be at least
$$m+h+(m+h)n +(m+h)n^2 + h + \bigl((m+h)n^2 +h\bigr)n  + B\bigl((m+h)n^2, n, h-1\bigr)$$ which may be rewritten as
$$(m+h)(n^3+n^2 +n +1) + h(n+1) + B\bigl((m+h)n^2, n, h-1\bigr).$$ Since this is larger than what was needed in Case (1), we no longer
need to consider Case (1).  This is the formula used in the recursive definition of $B(m,n,h).$\\

\noindent {\bf Subcase (2b):  $h'' = h$.}  In this case we show that $d = h$. To see this, note that the the relations on $\vect gn$ are
exactly given by $P$, and because $h'=h$ we know that the relations on $\vect \alpha n$ are exactly given by $P$  again.
But the elements $\vect \beta n$ are homomorphic images of the elements $\vect gn$,  so that these satisfy $P$, and so are the elements
$\alpha_1+\beta_1, \, \ldots, \alpha_n+\beta_n$, which must also satisfy $P$.  Moreover, since $h'' =h$, the algebra generated by
the elements $\vect \beta n$  has dimension $h$, and $P$ gives the ideal of all relations on the elements $\vect \beta n$.  We may
apply Lemma \ref{lin} to conclude that the $P$ is generated by linear forms. Since standard form was set up so that
$\vect f d$ are linearly independent but the $f_i$ for $i >h$  are integral over $K[\vect f h]$, we must have that
$f_i = 0$ for $i > h$,  i.e.,  $d = h$ as claimed.  Hence, the only nonzero tail polynomial are $\vect gh$, which form a regular
sequence modulo $\vect xm, \,\uu,\,\uv,\,\uw$.  Then $K[\vect F {m+n}] \inc K[\vect gh,\,\vect xm, \,\uu,\,\uv,\,\uw]$, and the
only constraint on $B(m,n,h)$ is that it be at least $(m+h)(1+n+n^2)$.  Therefore, using the formula obtained in
Subcase (2a) will cover all cases. $\square$ \\

\section{Size estimates}\label{size}

We study briefly the size of the functions $B(m,n,h)$, $C_0(s)$, and $C(s)$.
We shall use the alternative notation $B(m,n,h) = B_h(m,n) = B_h$.
Let $$g=g(n) = n^3+ n^2+n+1,$$
and let $$g_1=g_1(n) = g(n)+n+1.$$ Thus, $g$ and $g_1$ are both monic in $n$ of degree 3.
Let $\theta_h(t) = \sum_{i=0}^{h-1} t^i$, so that  $\theta_0 =0$  and $\theta_h = t\theta_{h-1}+1$ for $h \geq 1$.
We claim that for all $h$,   $$(m+h)\bigl((n+1)(n^2)^h + g\theta_h(n^2)\bigr)$$  is a lower bound for  $B_h$,
while $$(m+h)\bigl((n+1)(n^2+1)^{h} + g_1\theta_h(n^2+1)\bigr)$$ is an upper bound.  Note that the two agree with $B_0$ when  $h =0$.
The proof that the given expressions are bounds is by induction on $h$.  If we assume the estimates for a given  $h-1 \geq 0$,  we have\\
\quad\\
\noindent $B_h = (m+h)g + h(n+1) +B_{h-1}\bigl((m+h)n^2,n\bigr) $
$$\leq (m+h)g +(m+h)(n+1) + \bigl((m+h)n^2+h-1\bigr)\bigl((n+1)(n^2+1)^{h-1} + g_1 \theta_{h-1}(n^2+1)\bigr).$$
We may replace  the summand  $h-1$ in the first factor of the third term by  $m+h$.  The expression then becomes
$$(m+h)\biggl(g_1 + (n^2+1)\bigl((n+1)(n^2+1)^{h-1} + g_1 \theta_{h-1}(n^2+1)\biggr),$$ which simplifies to the required form.
We also have $$B_h \geq (m+h)g + h(n+1) +\bigl((m+h)n^2 +(h-1)\bigr)\bigl((n+1)(n^2)^{h-1} + g\theta_{h-1}(n^2)\bigr).$$
We may simply drop the terms  $h(n+1)$ and  $h-1$ to get the required estimate.\\

A lower bound for both $C_0(s)$ and for $C(s)$ may be obtained from the case $m = 0$, $s = n$, $h = s-1$, which yields
$$C(s) \geq (s-1)\bigl((s+1)(s^2)^{s-1} + g(s)\theta_{s-1}(s^2)\bigr) + (s-1).$$
Since whenever $h =s$ we must have $m=0$ and that the given forms are already a regular sequence, we
may assume $h \leq s-1$. This yields the upper bound   $$C(s) \leq (s-1)\bigl((s+1)(s^2+1)^{s-1} + g_1(s)\theta_{s-1}(s^2+1)\bigr) + (s-1).$$
Both bounds are asymptotic to  $2s^{2s}$,  since $\theta_{s-1}(t) = (t^{s-1}-1)/(t-1)$, and $(s^2+1)^{s-1}$ is asymptotic
to $(s^2)^{s-1}$. Hence, both $C_0(s)$ and $C(s)$ are asymptotic to $2s^{2s}$.

\section{The question in higher degree}
We raise here the following question:  is there an integer $C(n,d)$ such that given $n$ forms $F_1, \, \ldots, \, F_n$  of degree at most $d$ in a
polynomial ring $R$ over an infinite field $K$,
there exists a regular sequence of forms $G_1, \, \ldots, \, G_k \in R$ of degree at most $d$ with $k \leq C(n,d)$ such that
$F_1, \, \ldots, \, F_n \in K[G_1, \, \ldots, \, G_k]$?  Our main result here is, of course, the case $d=2$.  An affirmative answer would
yield  results analogous to our main theorems for arbitrary $d$:  $C(n,d)$ would bound the projective dimension of
$R/(F_1, \, \ldots, \, F_n)$ over any field, and $C(dn,d)$ would give corresponding results when the $F_i$ are polynomials of degree
at most $d$, not required to be homogeneous.
\bigskip

$\begin{array}{ll}
\textrm{Department of Mathematics}           &\quad \textrm{Department of Mathematics}\\
\textrm{University of Michigan}              &\quad \textrm{Adrian College}\\
\textrm{Ann Arbor, MI 48109--1109}           &\quad \textrm{Adrian, MI 49221}\\
\textrm{USA}                                 &\quad \textrm{USA}\\
\smallskip
\textrm{E-mail:}                             &\quad \textrm{E-mail:}\\
\textrm{hochster@umich.edu}                 &\quad \textrm{tananyan@adrian.edu}\\
\end{array}$

\end{document}